\documentclass[12pt,reqno]{amsart}
\usepackage{amssymb}
\usepackage{amsxtra}

\newcommand{\al}{\alpha}
\newcommand{\be}{\beta}
\newcommand{\C}{{\mathbf C}}
\newcommand{\ds}{\displaystyle}
\newcommand{\e}{\varepsilon}
\newcommand{\la}{\lambda}

\newcommand{\R}{{\mathbf R}}
\newcommand{\sigg}{\sigma_g}
\newcommand{\sigf}{\sigma_f}
\newcommand{\w}{\omega}
\newcommand{\Z}{{\mathbf Z}}
\newcommand{\z}{\zeta}
\newcommand{\eop}{\hfill$\square$}

\theoremstyle{plain}
\newtheorem{Thm}{Theorem}
\newtheorem{Prop}{Proposition}
\newtheorem{Cor}{Corollary}
\newtheorem{Lem}{Lemma}
\theoremstyle{definition}

\theoremstyle{remark}

\numberwithin{equation}{section}

\begin{document}

\title[Dirichlet Series Acceleration]{Series Acceleration Formulas for
Dirichlet Series with Periodic Coefficients}

\date{\today}

\author{David~M. Bradley}
\address{Department of Mathematics \& Statistics\\
         University of Maine\\
         5752 Neville Hall
         Orono, Maine 04469-5752\\
         U.S.A.}
\email{dbradley@e-math.ams.org, bradley@gauss.umemat.maine.edu}



\begin{abstract}
    Series acceleration formulas are obtained for
    Dirichlet series with periodic coefficients.  Special
    cases include Ramanujan's formula for the values of the Riemann
    zeta function at the odd positive integers exceeding two,
    and related formulas for values of Dirichlet $L$-series
    and the Lerch zeta function.
\end{abstract}

\maketitle

\section{Introduction}\label{sect:Intro}
Let $m$ be a positive integer and let $g$ be a complex-valued
function defined on the integers that is periodic with period $m$.
In other words, $g:\Z\to\C$ has the property that $g(n+m)=g(n)$
for all integers $n$. Examples include the constant functions
($m=1$) and the Dirichlet characters modulo $m>1$, but it is not
necessary in what follows to make any sort of multiplicativity
assumptions on $g$.  If $g:\Z\to\C$ has period $m$ then $g$ has
mean value
\begin{equation}
   M(g) = \frac1m \sum_{n=0}^{m-1} g(n),
\label{MeanValue}
\end{equation}
and furthermore, since $g$ is bounded ($\ds\sup_{n\in\Z}
|g(n)|=\max_{1\le n\le m} |g(n)|<\infty$), the Dirichlet series
\begin{equation}
   L(s,g) := \sum_{n=1}^\infty \frac{ g(n)}{n^s}
\label{Lsg}
\end{equation}
converges absolutely in the half-plane $\{s\in\C:\Re(s)>1\}$.  An
easy partial summation argument shows that the series~(\ref{Lsg})
converges conditionally in the  half-plane $\{s\in\C: \Re(s)>0\}$
if and only if $M(g)=0$.

Henceforth, we shall denote the abscissa of convergence of the
Dirichlet series~(\ref{Lsg}) by $\sigg$. That is,
\begin{equation}
   \sigg := \inf\{\sigma\in\R: (\ref{Lsg}) \mbox{ converges for }
   \Re(s)>\sigma\}.
\label{siggdef}
\end{equation}
Then the previous observations may be restated as asserting that
$\sigg\le 1$ for all periodic $g:\Z\to\C$, and $\sigg\le 0$ if and
only if $M(g)=0$.  Of course if $g:\Z\to\C$ is periodic and
$M(g)=0$, then we actually must have $\sigma_g=0$ except in the
trivial case when $g$ vanishes identically and $\sigma_g=-\infty$.

As part of a general program aimed at developing methods for
calculating important number-theoretical constants to high
precision, we consider here the problem of replacing values of the
series~(\ref{Lsg}) by equivalent expressions with improved rate of
convergence.  If we succeed in finding such an expression for a
specific value of $s$, we say we have obtained a {\em series
acceleration formula} for $L(s,g)$.  Theorems 1 and 2 below
provide series acceleration formulas for $L(s,g)$ when $g:\Z\to\C$
is periodic and either odd or even and $s$ is a positive integer.
Special cases include Ramanujan's beautiful reciprocity formula
for the values $\z(2q+1)$ of the Riemann zeta function (where $q$
is a positive integer), analogous formulas for values of Dirichlet
$L$-series, and other related results that have appeared in the
literature.

We make no attempt here to give a rigorous definition of the
concept of a series acceleration formula. However, the following
remarks should give a reasonable indication of what we mean.  All
our results have the form $\sum a_n = C + S$, where $C$ is a
``closed form'' expression and $S$ is a finite sum of series
(equivalently, a single series) of the form $\sum b_n$, in which
$b_n$ is an elementary function of $n$ and $ \limsup_{n\to\infty}
|b_n/a_n|^{1/n} <1$.  Of course, ``closed form'' depends on what
sort of objects one is prepared to accept as fundamental.  Let us
agree to accept values of the elementary functions at the integers
as closed form.  Then (see Proposition~\ref{Prop:Lclosedcplx})
$L(q,g)$ is closed form when $q$ is a positive integer such that
$g$ and $q$ are of the same parity, and in particular, the numbers
$\z(2q)$ are closed form.  For even $g$, the second formula in
Theorem~\ref{Thm:gcplx} expresses $L(2q+1,g)$ in terms of the
aforementioned closed form values, rapidly convergent series, and
$\z(2q+1)$. Since we do not regard $\z(2q+1)$ as closed form, a
series acceleration formula for $\z(2q+1)$ is needed in order to
achieve a series acceleration formula for $L(2q+1,g)$.
Fortunately, Ramanujan's formula for $\z(2q+1)$ (see
Corollary~\ref{Cor:Z}) serves quite adequately as a series
acceleration formula, and thus Theorems~\ref{Thm:gcplx}
and~\ref{Thm:greal} do indeed provide legitimate series
acceleration formulas for $L(s,g)$, as claimed.

\section{Main Result}

\begin{Thm}\label{Thm:gcplx}
Let $m$ and $q$ be positive integers, $\w=\exp(2\pi i/m)$,
$g:\Z\to\C$ periodic of period $m$, $M(g)$ as in~(\ref{MeanValue})
and $L(s,g)$ as in~(\ref{Lsg}). If $\al$ and $\be$ are positive
real numbers satisfying $\al\be=\pi^2$, then
\begin{eqnarray*}
   \al^{-q+1/2}\bigg\{\tfrac12 L(2q,g) + \sum_{n=1}^\infty
   \frac{n^{-2q}g(n)}{e^{2n\al}-1}\bigg\}
   &=& (-1)^q \be^{-q+1/2} i m^{-1} \sum_{k=1}^{m-1}
   g(k)\sum_{n=1}^\infty
   \frac{n^{-2q}}{e^{2n\be/m}\w^k -1}\\
   &+& \sum_{j=0}^q
   (-1)^{j+1}\al^{j-1/2-q}\be^{-j}\z(2j)L(2q-2j+1,g)
\end{eqnarray*}
if $g$ is odd, and
\begin{multline*}
   \al^{-q}\bigg\{\tfrac12 L(2q+1,g)+\sum_{n=1}^{\infty}
   \frac{n^{-2q-1}g(n)}{e^{2n\al}-1}\bigg\}\\
   =  (-\be)^{-q}\bigg\{\tfrac12 M(g)\z(2q+1)+
      \frac1m\sum_{k=0}^{m-1}g(k)
      \sum_{n=1}^\infty \frac{n^{-2q-1}}{e^{2n\be/m}\w^k-1}\bigg\}\\
   + \sum_{j=0}^{q+1}
   (-1)^{j+1}\al^{j-q-1}\be^{-j}\z(2j)L(2q-2j+2,g)
\end{multline*}
if $g$ is even. If $M(g)=0$, the latter formula is also valid when
$q=0$.
\end{Thm}

There is an equivalent version of Theorem~\ref{Thm:gcplx} in which
the terms involving complex $m$th roots of unity are paired so as
to yield a real-valued expression:

\begin{Thm}\label{Thm:greal}
Let $m$ and $q$ be positive integers, $g:\Z\to\C$ periodic of
period $m$, $M(g)$ as in~(\ref{MeanValue}) and $L(s,g)$ as
in~(\ref{Lsg}). If $\al$ and $\be$ are positive real numbers
satisfying $\al\be=\pi^2$, then
\begin{multline*}
   \al^{-q+1/2}\bigg\{\tfrac12 L(2q,g) + \sum_{n=1}^\infty
   \frac{n^{-2q}g(n)}{e^{2n\al}-1}\bigg\}\\
   = \tfrac12 (-1)^q \be^{-q+1/2} m^{-1} \sum_{k=1}^{m-1}
   g(k)\sin(2\pi k/m) \sum_{n=1}^\infty
   \frac{n^{-2q}}{\cosh(2n\be/m)-\cos(2\pi k/m)}\\
   + \sum_{j=0}^q
   (-1)^{j+1}\al^{j-1/2-q}\be^{-j}\z(2j)L(2q-2j+1,g)
\end{multline*}
if $g$ is odd, and
\begin{multline*}
   \al^{-q}\bigg\{\tfrac12 L(2q+1,g)+\sum_{n=1}^{\infty}
   \frac{n^{-2q-1}g(n)}{e^{2n\al}-1}\bigg\}\\
   = \tfrac12 (-\be)^{-q} m^{-1}\sum_{k=0}^{m-1}g(k)
      \sum_{n=1}^\infty \bigg(\frac{\cos(2\pi k/m)-\exp(-2n\be/m)}
      {\cosh(2n\be/m)-\cos(2\pi k/m)}\bigg)n^{-2q-1}\\
   +\tfrac12
   (-\be)^{-q}M(g)\z(2q+1)
   + \sum_{j=0}^{q+1}
   (-1)^{j+1}\al^{j-q-1}\be^{-j}\z(2j)L(2q-2j+2,g)
\end{multline*}
if $g$ is even.  If $M(g)=0$, the latter formula is also valid
when $q=0$.
\end{Thm}

Our proof of Theorems~\ref{Thm:gcplx} and~\ref{Thm:greal} is
outlined in Section~\ref{sect:ProveMain}.  It should be noted that
Ramanujan had a generalization of Theorem~\ref{Thm:gcplx} in which
$g$ is replaced by an entire function satisfying suitable growth
conditions.  Ramanujan's generalization is proved using contour
integration in~\cite[pp.\ 429--430]{BerndtIV}.  Although our
results are less general than Ramanujan's, our method of proof is
somewhat more elementary than previous approaches, and therefore
may be of some interest.

As we remarked in the Introduction, it should be noted that for
these results to be of use as series acceleration formulas, it is
necessary to have a series acceleration formula for $\z(2q+1)$ in
the case when $g$ is even and $M(g)\ne 0$, and a closed-form
evaluation of $L(s,g)$ when $s$ is a non-negative integer of the
same parity as $g$.  For $\z(2q+1)$, we may use Ramanujan's
formula---see Corollary~\ref{Cor:Z} below. For $L(s,g)$, we first
recall that if $g=\chi$ is a Dirichlet character and $s$ is a
positive integer such that $(-1)^s\chi(-1)=1$, then the
corresponding Dirichlet $L$-series $L(s,\chi)$ has a closed-form
evaluation in terms of the Gauss sum
\[
   G(\chi) = \sum_{k=0}^{m-1} \chi(k) e^{2\pi i k/m}
\]
and the so-called {\em generalized Bernoulli numbers}
$B_{n,\chi}$, defined for non-negative integers $n$ by the formula
\[
   \sum_{k=1}^m \chi(k)\frac{te^{kt}}{e^{mt}-1}
   = \sum_{n=0}^\infty  B_{n,\chi}\,\frac{t^n}{n!}.
\]
This well-known result is typically derived (see eg.~\cite{Jurg})
as a consequence of the functional equation relating $L(s,\chi)$
to $L(1-s,\overline{\chi})$, but actually one does not need the
functional equation to accomplish this, and in fact, one does not
even need the multiplicativity property of the Dirichlet
characters.

As customary, define the Bernoulli polynomials $B_n(x)$ by their
exponential generating function
\begin{equation}
   \frac{te^{xt}}{e^t-1} = \sum_{n=0}^\infty B_n(x)\frac{t^n}{n!},
   \qquad |t|<2\pi,
\end{equation}
and the Bernoulli numbers by $B_n= B_n(0)$ for $0\le n\in\Z$. The
following closed-form evaluation for $L(s,g)$ makes no
multiplicativity assumption on the periodic arithmetical function
$g$.

\begin{Prop}[\cite{BerndtSchoen}]\label{Prop:Lclosedcplx} Let $m$ be a positive
integer, $g:\Z\to\C$ periodic with period $m$ and let $L(s,g)$ be
as in~(\ref{Lsg}).  If $q$ is a non-negative integer such that $g$
and $q$ are both odd or both even, then
\[
   L(q,g) = -\frac12\cdot
   \frac{(2\pi i)^q}{q!}\sum_{k=0}^{m-1} \widehat{g}(k)
   B_q(k/m),
\]
where $\widehat{g}(k)$ is the $k$th discrete Fourier coefficient
of $g$ defined by
\[
   \widehat{g}(k):= \frac{1}{m}\sum_{j=0}^{m-1}g(j)\, e^{-2\pi
   ijk/m}.
\]
\end{Prop}

\noindent{\bf Proof.}  See equations (6.23) and (6.25)
of~\cite{BerndtSchoen} for the case when $q$ is a positive
integer.  For the case $q=0$, we must show that if $g$ is even,
then
\begin{eqnarray*}
   L(0,g) &=& -\frac12 \sum_{k=0}^{m-1} \widehat{g}(k)
   = -\frac{1}{2m}\sum_{k=0}^{m-1} \sum_{j=0}^{m-1} g(j) e^{-2\pi
   ijk/m}
   = -\frac{1}{2m}\sum_{j=0}^{m-1}g(j)\sum_{k=0}^{m-1} e^{-2\pi
   ijk/m}\\
   &=& -\tfrac12 g(0).
\end{eqnarray*}
But this one can easily establish by meromorphically extending the
definition of $L(s,g)$ to the half-plane $\{s\in\C: \Re(s)>-1\}$.
\eop

The equivalent ``real'' version of
Proposition~\ref{Prop:Lclosedcplx} is stated below for
convenience.

\begin{Cor}\label{Lclosedreal} Let $m$ be a positive
integer, $g:\Z\to\C$ periodic with period $m$ and let $L(s,g)$ be
as in~(\ref{Lsg}).   Then for all non-negative integers $q$,
\[
   L(2q,g) = \frac{(-1)^{q+1}(2\pi)^{2q}}{(2q)!}\cdot
   \frac{1}{2m}\sum_{k=0}^{m-1}B_{2q}(k/m)\sum_{j=0}^{m-1}
   g(j)\cos(2\pi j k/m)
\]
if $g$ is even, and
\[
   L(2q+1,g) = \frac{(-1)^{q+1} (2\pi)^{2q+1}}{(2q+1)!}\cdot
   \frac{1}{2m}\sum_{k=0}^{m-1}B_{2q+1}(k/m)\sum_{j=0}^{m-1}
   g(j)\sin(2\pi jk/m)
\]
if $g$ is odd.
\end{Cor}

In particular, we have the following well-known evaluations which
are needed in the sequel:

\begin{Cor}
For $\Re(s)>0$, let $L(s):=\sum_{k=0}^\infty (-1)^k (2k+1)^{-s}$,
and let $\z(s)$ denote the Riemann zeta function. If $n$ is a
non-negative integer, then
\begin{equation}
   \z(2n) = -\frac12\cdot\frac{(2\pi i)^{2n}B_{2n}}{(2n)!}
\label{z2n}
\end{equation}
and
\begin{equation}
   L(2n+1) = \frac12\left(\frac{\pi}{2}\right)^{2n+1}\frac{(-1)^{n} E_{2n}}{(2n)!},
\label{Lodd}
\end{equation}
where
\[
   \frac{1}{\cosh(t)} = \sum_{k=0}^\infty \frac{E_k}{k!}t^k,
   \qquad |t|<\tfrac12\pi,
\]
generates the Euler numbers $E_k$ for non-negative integers $k$.
\end{Cor}

\section{Consequences}\label{sect:Consequences}

\begin{Cor}[Ramanujan's formula for $\z(2q+1)$]
\label{Cor:Z} Let $q$ be a positive integer, and let
$\al$ and $\be$ be positive real numbers with $\al\be=\pi^2$.
Then
\begin{multline}
\label{RamaZodd}
   \al^{-q}\bigg\{\tfrac12\z(2q+1)+\sum_{n=1}^\infty
   \frac{n^{-2q-1}}{e^{2n\al}-1}\bigg\}
   = (-\be)^{-q}\bigg\{\tfrac12\z(2q+1)+\sum_{n=1}^\infty
   \frac{n^{-2q-1}}{e^{2n\be}-1}\bigg\}\\
   + 2^{2q}\sum_{k=0}^{q+1}(-1)^{k+1}\frac{B_{2k}}{(2k)!}
   \frac{B_{2q+2-2k}}{(2q+2-2k)!}\,\al^{q+1-k}\be^{k}.
\end{multline}
\end{Cor}

\noindent{\bf Proof.}  Let $m=1$ in Theorem~\ref{Thm:gcplx} and
let $g(0)=1$, noting that periodicity implies $g\equiv 1$ is
constant. The stated formula now follows after appropriately
substituting~(\ref{z2n}) and performing standard algebraic
manipulations. \eop

\noindent{\bf Remarks.} Corollary~\ref{Cor:Z} appears in
Ramanujan's notebooks~\cite[vol. I, p.\ 259, no.\ 15; vol. II, p.\
177, no.\ 21] {RamaBooks}, but he did not publish a proof.
Although Lerch~\cite{Lerch} proved the special case $\al=\be=\pi$,
Grosswald's paper~\cite{Gross} was responsible for generating much
of the subsequent interest in~(\ref{RamaZodd}).
Guinand~\cite{Guin} showed how Corollary~\ref{Cor:Z} arises from
the modular transformation $V(z)=-1/z$ of the function
\[
   f(z,-2n)  = \sum_{k=1}^\infty \frac{k^{-2n-1}}{e^{-2\pi
   ikz}-1}, \qquad \Im(z)>0,
\]
and letting $z=i$.  Additional formulas for $\z(2q+1)$ and cognate
results for certain other Dirichlet series may be obtained by
applying other transformations. See~\cite{BerndtRocky} for a
comprehensive account with extensive references to the literature.
Many further references can be found in~\cite[p.\ 276]{BerndtII}.
More recently, Ramanujan's formula for $\z(2q+1)$ has been
applied~\cite{Helmut} in studying the variance of the random
variable $X_q$ representing the number of internal nodes of a
binary trie built from $q$ data.

\begin{Cor}\label{Cor:L4}
For $\Re(s)>0$, let $L(s):=\sum_{n=0}^\infty (-1)^n (2n+1)^{-s}$.
Let $q$ be a positive integer, and let $\al$ and $\be$ be positive
real numbers with $\al\be=\pi^2$.  Then
\begin{multline*}
   \al^{-q+1/2}\bigg\{\tfrac12
   L(2q)+\sum_{n=0}^{\infty}\frac{(-1)^n
   (2n+1)^{-2q}}{e^{(4n+2)\al}-1}\bigg\}\\
   = \tfrac14 (-1)^q \be^{-q+1/2}\sum_{n=1}^\infty \frac{n^{-2q}}
   {\cosh(n\be/2)}
   +2^{2q-3}\sum_{k=0}^q (-1)^k 2^{-4k} \frac{E_{2k}}{(2k)!}
    \frac{B_{2q-2k}}{(2q-2k)!}\al^{q-k}\be^{k+1/2}.
\end{multline*}
\end{Cor}

\noindent{\bf Proof.} Put $m=4$ in Theorem~\ref{Thm:greal}  and
let $g(2n+1)=(-1)^n$ and $g(2n)=0$ for $n\in\Z$.  The stated
formula now follows after appropriately substituting~(\ref{z2n})
and~(\ref{Lodd}) and performing standard algebraic
manipulations.\eop

\noindent{\bf Remarks.} Corollary~\ref{Cor:L4} corrects the
misprints in Proposition 3.5 of~\cite[p.\ 169]{BerndtRocky}.  As
noted by Berndt~\cite{BerndtRocky}, the result was known to
Ramanujan~\cite[Vol.\ I, p.\ 274; Vol.\ II, pp.\
177--178]{RamaBooks}, but the first published proof is due to
Chowla~\cite{Chowla}.  Note that $L(s)$ is the Dirichlet
$L$-series corresponding to the primitive non-principal Dirichlet
character modulo 4.  In a similar vein, formulas for general
Dirichlet $L$-functions have been given by Katayama~\cite{Kat} and
Berndt~\cite{BerndtActa}. Of course, since all Dirichlet
$L$-series have periodic coefficients, our
Theorems~\ref{Thm:gcplx} and~\ref{Thm:greal} can be specialized to
give ``Ramanujan formulas'' for general Dirichlet $L$-functions as
well.

For purposes of maximizing the rate of convergence of the most
slowly convergent series in Theorems~\ref{Thm:gcplx}
and~\ref{Thm:greal}, the optimal choice of $\al$ and $\be$ is
$\al=\pi/\sqrt{m}$, $\be=\pi\sqrt{m}$, respectively.  If this
choice is made in Corollary~\ref{Cor:L4} with $\sqrt{m}=q=2$, we
recover Ramanujan's formula~\cite[p.\ 43]{Ti2} for Catalan's
constant:
\begin{equation}
   G := \sum_{n=0}^\infty \frac{(-1)^n}{(2n+1)^2}
   = \frac{5}{48}\pi^2 - 2\sum_{n=0}^\infty \frac{(-1)^n
   (2n+1)^{-2}}{e^{(2n+1)\pi}-1}-\frac{1}{4}\sum_{n=1}^\infty
   \frac{1}{n^2\cosh(\pi n)}.
\label{Catalan}
\end{equation}

\begin{Cor}[Theorem 3.3 of~\cite{BerndtRocky}]
\label{Cor:sincos} Let $q$ be a non-negative integer and let
$\al$, $\be$ and $r$ be positive real numbers with $0<r<1$ and
$\al\be=\pi^2$.   Then
\begin{multline}
\label{cos}
   \al^{-q}\bigg\{\frac12\sum_{n=1}^\infty n^{-2q-1}\cos(2\pi n
   r)+\sum_{n=1}^\infty \frac{n^{-2q-1}\cos(2\pi
   nr)}{e^{2n\al}-1}\bigg\}\\
   = (-\be)^{-q}\bigg\{\frac12\sum_{n=1}^\infty
   n^{-2q-1}e^{-2n\be r}+\sum_{n=1}^\infty
   \frac{n^{-2q-1}\cosh(2 n\be r)}{e^{2n\be}-1}\bigg\}\\
   -2^{2q}\sum_{k=0}^{q+1}(-1)^k \frac{B_{2k}(r)}{(2k)!}
   \frac{B_{2q+2-2k}}{(2q+2-2k)!}\,\al^{q+1-k}\be^{k};
\end{multline}
and if $q$ is a positive integer, then
\begin{multline}
\label{sin}
   \al^{-q+1/2}\bigg\{\frac12\sum_{n=1}^\infty n^{-2q}\sin(2\pi
   nr)+\sum_{n=1}^\infty \frac{n^{-2q}\sinh(2\pi
   nr)}{e^{2n\al}-1}\bigg\}\\
   =(-1)^q\be^{-q+1/2}\bigg\{\frac12 \sum_{n=1}^\infty n^{-2q}e^{-2n\be r}
   +\sum_{n=1}^\infty \frac{n^{-2q}\sinh(2n\be
   r)}{e^{2n\be}-1}\bigg\}\\
   -2^{2q-1}\sum_{k=0}^q (-1)^k\frac{B_{2k+1}(r)}{(2k+1)!}
   \frac{B_{2q-2k}}{(2q-2k)!}\,\al^{q-k}\be^{k+1/2}.
\end{multline}
\end{Cor}

\noindent{\bf Proof.}  It suffices to prove the given formulas in
the case when $0<r<1$ and $r$ is rational, as the result for real
$r$ with $0<r<1$ then follows by taking limits.   We shall prove
only~(\ref{cos}), as the proof of~(\ref{sin}) is almost identical.
Suppose $u$ and $m$ are integers satisfying $0<u<m$. Let $r=u/m$,
and let $g:\Z\to\C$ be defined by $g(n)=\cos(2\pi nr)$ for all
integers $n$.  Then $g$ is even and periodic with period $m$.
Rewrite Theorem~\ref{Thm:gcplx} in the form
\begin{multline*}
   \tfrac12 \al^{-q} \sum_{n=1}^\infty n^{-2q-1}\coth(n\al)
   = \tfrac 12(-\be)^{-q}
   \sum_{n=1}^\infty n^{-2q-1}\coth((n\be + \pi i k)/m)\\
   - \sum_{k=0}^{q+1} (-1)^{q+1-k}
   L(2k,g)\z(2q-2k+2)\al^{-k} \be^{k-q-1}.
\end{multline*} From Proposition~\ref{Prop:Lclosedcplx}, we have
\[
   L(2k,g) = -\frac12\cdot \frac{(2\pi i)^{2k} B_{2k}(r)}{(2k)!},
   \qquad 0\le k\in\Z.
\]
Using also~(\ref{z2n}), we see that it now suffices (with
$\la=n\be$) to prove that for all $\la>0$,
\[
   m^{-1}\sum_{k=0}^{m-1}\cos(2\pi k r) \coth((\la +\pi i k)/m)
   = e^{-2\la r} + \frac{2\cosh(2\la r)}{e^{2\la}-1}
   = \frac{\cosh((1-2r)\la)}{\sinh(\la)}.
\]
But
\begin{eqnarray*}
  && m^{-1}\sum_{k=0}^{m-1}\cos(2\pi k r) \coth((\la +\pi i
  k)/m)\\
   &=& \frac{1}{2m}\sum_{k=0}^{m-1}\big(e^{2\pi i ku/m}+e^{-2\pi
   iku/m}\big)\frac{1+e^{-2(\la +\pi ik)/m}}{1-e^{-2(\la +\pi
   ik)/m}}\\
   &=& \frac{1}{2m}\sum_{k=0}^{m-1}\big(e^{2\pi i ku/m}+e^{-2\pi
   iku/m}\big)\big(1+e^{-2(\la +\pi ik)/m}\big)
   \sum_{j=0}^\infty e^{-2(\la+\pi ik)j/m}\\
   &=&  \frac{1}{2m}\sum_{j=0}^\infty e^{-2\la j/m}
   \sum_{k=0}^{m-1}\big(e^{2\pi i k(u-j)/m}+e^{-2\pi
   ik(u+j)/m}\big)\\
   &+& \frac{1}{2m}\sum_{j=1}^\infty e^{-2\la j/m}
   \sum_{k=0}^{m-1}\big(e^{2\pi i k(u-j)/m}+e^{-2\pi
   ik(u+j)/m}\big)\\
   &=& \sum_{k=0}^\infty \big(
   e^{-2\la(u+km)/m}+e^{-2\la(m-u+km)/m}\big)\\
   &=& \frac{e^{-2\la r}}{1-e^{-2\la}}+\frac{e^{-2\la(1-r)}}{1-e^{-2\la}}
   = \frac{e^{-(1-2r)\la}+e^{(1-2r)\la}}{e^{\la}-e^{-\la}}\\
   &=& \frac{\cosh((1-2r)\la)}{\sinh(\la)}.
\end{eqnarray*}
This completes the proof of~(\ref{cos}).  The proof of~(\ref{sin})
proceeds mutatis mutandis with $g(n)=\sin(2\pi nr)$ replacing
$g(n)=\cos(2\pi nr)$.  \eop

\noindent{\bf Remarks.} In~(\ref{sin}), we have corrected the
misprints in the corresponding formula (3.12) of~\cite[p.\
167]{BerndtRocky}.  Berndt~\cite{BerndtPeriodic} has given a
generalization of Corollary~\ref{Cor:sincos} to periodic
sequences.  In~\cite{BerndtRocky}, Berndt deduces several
interesting results by specializing~(\ref{cos}) and~(\ref{sin}) in
various ways. Here, we confine ourselves to remarking that setting
$r=1/4$ gives Corollary~\ref{Cor:L4}, and letting $r$ tend to zero
gives Euler's formula~(\ref{z2n}) and Ramanujan's
formula~(\ref{RamaZodd}) again.

\section{Proof of Main Result}\label{sect:ProveMain}

Although the results of section~\ref{sect:Consequences} have been
previously derived as consequences of quite general modular
transformation formulas~\cite{BerndtRocky}, we feel that it may
nevertheless be of interest to give a proof of our main result
using elementary methods analogous to Ramanujan's~\cite{Ti2} proof
of~(\ref{Catalan}).  The main idea is to employ the partial
fraction expansion of the hyperbolic cotangent in a non-trivial
manner.  The proof has been broken down into a sequence of
relatively straightforward lemmata, the last of which gives our
main result after making a trivial substitution.

\begin{Lem}\label{Lem:1}
Let $f:\Z^+\to\C$ and suppose that the associated Dirichlet series
$F(s):=\sum_{n=1}^\infty n^{-s}f(n)$ has abscissa of convergence
$\sigf<\infty$.  If $\Re(s)>\sigf-1$, $g:\Z\to\C$ and $x$ is any
positive real number, then
\[
   \tfrac12 \pi x F(s+1) + \pi x\sum_{k=1}^\infty
   \frac{k^{-s-1}f(k)}{e^{2\pi k/x}-1}
   = \tfrac12 x^2 F(s+2)+ x^2\sum_{n=1}^\infty\sum_{k=1}^\infty
   \frac{k^{-s}f(k)}{k^2+n^2x^2}.
\]
\end{Lem}

\noindent{\bf Proof.}  Recalling the partial fraction expansion of
the hyperbolic cotangent~\cite[p.\ 259, 6.3.13]{AS}, we find that
\begin{eqnarray*}
   \tfrac12\pi x F(s+1) + \pi x\sum_{k=1}^\infty
      \frac{k^{-s-1}f(k)}{e^{2\pi k/x}-1}
   &=& \tfrac12\pi x\sum_{k=1}^\infty k^{-s-1}f(k)\bigg(1+
      \frac{2}{e^{2\pi k/x}-1}\bigg)\\
   &=&\tfrac12\pi x\sum_{k=1}^\infty k^{-s-1} f(k)\coth(\pi k/x)\\
   &=&\sum_{k=1}^\infty k^{-s}f(k)\bigg(\tfrac12 x^2 k^{-2}
       +\sum_{n=1}^\infty \frac{x^2}{k^2+n^2x^2}\bigg)\\
   &=&\tfrac12 x^2 F(s+2)+\sum_{k=1}^\infty k^{-s}f(k)
       \sum_{n=1}^\infty \frac{x^2}{k^2+n^2x^2}\\
   &=&\tfrac12 x^2 F(s+2)+ x^2\sum_{n=1}^\infty
   \sum_{k=1}^\infty \frac{k^{-s}f(k)}{k^2+n^2x^2}.
\end{eqnarray*}
The interchange of summation in the final step can be justified by
the Fubini-Tonelli theorem with counting measure~\cite[p.\
67]{Folland}.  The double sum is absolutely convergent because if
$0<\e<\min(1,\Re(s)+1-\sigf)$, then concavity of the logarithm
implies that
\[
   \frac{1}{k^2+n^2x^2}
   \le \frac{2}{(1-\e)k^2+(1+\e)n^2x^2}
   \le \frac{1}{k^{1-\e}(nx)^{1+\e}},
\]
and so
\[
   \sum_{n,k=1}^\infty \bigg|\frac{k^{-s}f(k)}{k^2+n^2x^2}\bigg|
   \le \sum_{n=1}^\infty \frac{1}{|nx|^{1+\e}}
   \sum_{k=1}^\infty \bigg|\frac{f(k)}{k^{s+1-\e}}\bigg|<\infty,
\]
since $\Re(s)+1-\e>\sigf$. \eop

\begin{Lem}\label{Lem:2}
Let $f:\Z^+\to\C$ and suppose that the associated Dirichlet series
$F(s):=\sum_{n=1}^\infty n^{-s} f(n)$ has abscissa of convergence
$\sigf<\infty$.  For $\Re(s)>\sigf-2$ and $y$ real, define
\begin{equation}
   T_f(s,y) := \sum_{k=1}^\infty \frac{k^{-s}f(k)}{k^2+y^2}.
\label{Tdef}
\end{equation}
If $q$ is a non-negative integer, $\Re(s)>\sigf+2q-2$ and $y\ne
0$, then
\[
   T_f(s,y) = (-1)^q y^{-2q}T_f(s-2q,y)+\sum_{j=1}^q (-1)^{j+1}
   y^{-2j}F(s-2j+2).
\]
\end{Lem}

\noindent{\bf Proof.} First, note that $T_f(s,0)=F(s+2)$, so the
series~(\ref{Tdef}) defining $T_f(s,y)$ converges if
$\Re(s)>\sigf-2.$  Next, observe that the formula clearly holds if
$q=0$, so we may assume that $q$ is a positive integer.  If
$\Re(s)>\sigf$ and $y\ne 0$, then
\begin{equation}
   T_f(s,y) = y^{-2}\sum_{k=1}^\infty
   \left(k^{-2}-(k^2+y^2)^{-1}\right)k^{-s+2}f(k)
   = y^{-2}F(s)-y^{-2}T_f(s-2,y).
\label{Trecurrence}
\end{equation}
The stated formula now follows by iterating~(\ref{Trecurrence}) or
alternatively by replacing $s$ by $s-2j+2$, multiplying both sides
by $(-1)^{j+1}y^{-2j+2}$ and telescoping the sum on $j$ from $1$
to $q$. \eop

\begin{Lem}\label{Lem:3}
Let $f:\Z^+\to\C$ and suppose that the associated Dirichlet series
$F(s):=\sum_{n=1}^\infty n^{-s}f(n)$ has abscissa of convergence
$\sigf<\infty$.  Let $q$ be a positive integer, $x$ a positive
real number, and $T_f$ as in~(\ref{Tdef}). Then for
$\Re(s)>\sigf+2q-2$, we have
\begin{multline*}
   \tfrac12\pi x F(s+1)+\pi x\sum_{n=1}^\infty
   \frac{n^{-s-1}f(n)}{e^{2\pi n/x}-1}\\
   = \sum_{j=0}^q (-1)^{j+1}x^{-2j+2}\z(2j)F(s-2j+2)
   + (-1)^q x^{-2q+2}\sum_{n=1}^\infty n^{-2q}T_f(s-2q,nx).
\end{multline*}
\end{Lem}

\noindent{\bf Proof.} By Lemma~\ref{Lem:1} and Lemma~\ref{Lem:2},
we have
\begin{eqnarray*}
   &&\tfrac12\pi x F(s+1)+\pi x\sum_{n=1}^\infty
   \frac{n^{-s-1}f(n)}{e^{2\pi n/x}-1}\\
   &=& \tfrac12 x^2 F(s+2)+x^2\sum_{n=1}^\infty T_f(s,nx)\\
   &=& \tfrac12 x^2 F(s+2)+x^2\sum_{n=1}^\infty \bigg\{(-1)^q
   n^{-2q}x^{-2q}\,T_f(s-2q,nx)\\
   && +\sum_{j=1}^q(-1)^{j+1}n^{-2j}x^{-2j}
   F(s-2j+2)\bigg\}\\
   &=&  \tfrac12 x^2 F(s+2) + \sum_{j=1}^q (-1)^{j+1}\z(2j)
      x^{-2j+2}F(s-2j+2)\\
   && +(-1)^qx^{-2q+2}\sum_{n=1}^\infty
      n^{-2q}\, T_f(s-2q,nx).
\end{eqnarray*}
Since $\z(0)=-1/2$, the sum on $j$ can be extended to absorb the
term $\tfrac12 x^2 F(s+2)$, and the stated formula follows. \eop

\begin{Lem}\label{Lem:4}
Let $m$ be a positive integer, $y$ a positive real number, and
$g:\Z\to\C$ periodic of period $m$.  Then
\[
   P.V. \sum_{0\ne k\in\Z}\frac{g(k)}{y+ik}
   = -y^{-1}g(0)+\pi m^{-1}\sum_{k=0}^{m-1}
   g(k)\coth(\pi(y+ik)/m),
\]
where ``P.V.'' denotes the principal value, i.e.\ the symmetric
limit
\[
   \lim_{N\to\infty} \sum_{0<|k|<N}\frac{g(k)}{y+ik}.
\]
\end{Lem}

\noindent{\bf Proof.}  Since $g$ has period $m$, we have
\begin{eqnarray*}
   P.V.\sum_{0\ne k\in\Z}\frac{g(k)}{y+ik}
   &=& -y^{-1}g(0)+\sum_{k=0}^{m-1}
   P.V.\sum_{r\in\Z}\frac{g(k+mr)}{y+i(k+mr)}\\
   &=& -y^{-1}g(0)+\sum_{k=0}^{m-1}g(k)\,
      P.V. \sum_{r\in\Z}\frac{1/m}{(y+ik)/m+ir}\\
   &=& -y^{-1}g(0)+\pi m^{-1}\sum_{k=0}^{m-1}g(k)\coth(\pi(y+ik)/m),\\
\end{eqnarray*}
as stated. \eop

\begin{Lem}\label{Lem:5}
Let $m$ be a positive integer, $\w:=\exp(2\pi i/m)$, $y$ a
positive real number, $g:\Z\to\C$ odd and periodic of period $m$,
$L(s,g)$ as in~(\ref{Lsg}), and $\sigg$ as in~(\ref{siggdef}).  As
in~(\ref{Tdef}), for $\Re(s)>\sigg-2$, define
\begin{equation}
   T_g(s,y) := \sum_{k=1}^\infty \frac{k^{-s}g(k)}{k^2+y^2}.
\label{Tgdef}
\end{equation}
Then
\begin{eqnarray}
   T_g(-1,y) = \sum_{k=1}^\infty \frac{kg(k)}{k^2+y^2}
   &=& \tfrac12\pi i m^{-1} \sum_{k=1}^{m-1} g(k)\coth(\pi(y+ik)/m)
   \label{ToddCoth}\\
   &=& \frac{\pi i}{m} \sum_{k=1}^{m-1} \frac{g(k)}{e^{2\pi y/m}
       \w^k-1}\label{ToddOmega}\\
   &=& \frac{\pi}{2m}\sum_{k=1}^{m-1}
     \frac{g(k)\sin(2\pi k/m)}{\cosh(2\pi y/m)-\cos(2\pi k/m)}
     \label{ToddReal}.
\end{eqnarray}
\end{Lem}

\noindent{\bf Proof.} We may as well assume $g$ does not vanish
identically, for otherwise the result is completely trivial. Since
$g$ has period $m$, $g(0)=g(m)=g(-m)$. Thus, as $g$ is odd,
$g(0)=\tfrac12(g(m)+g(-m))=0$.  Furthermore,
\[
   2mM(g) =
   \sum_{k=0}^{m-1}(g(k)+g(m-k))=\sum_{k=0}^{m-1}(g(k)+g(-k))=0.
\]
Therefore, by the remarks in the second paragraph of the
Introduction, $\sigg=0$, and the series defining $T_g(-1,y)$
converges.  But since $g$ is odd,
\begin{eqnarray*}
   \sum_{k=1}^\infty \frac{k g(k)}{k^2+y^2}
   &=& \frac{i}{2}\sum_{k=1}^\infty \bigg(\frac{1}{y+ik}
      -\frac{1}{y-ik}\bigg)g(k)\\
   &=& \frac{i}{2}\sum_{k=1}^\infty \bigg(\frac{g(k)}{y+ik}+
       \frac{g(-k)}{y-ik}\bigg)\\
   &=& \frac{i}{2}P.V.\sum_{0\ne k\in\Z} \frac{g(k)}{y+ik}\\
   &=& \tfrac12\pi i m^{-1} \sum_{k=1}^{m-1} g(k)\coth(\pi(y+ik)/m)
\end{eqnarray*}
by Lemma~\ref{Lem:4} and the fact that $g(0)=0$.  This
establishes~(\ref{ToddCoth}).

Next, since $g(0)=M(g)=0$,
\begin{eqnarray*}
   \tfrac12\pi i m^{-1} \sum_{k=1}^{m-1} g(k)\coth(\pi(y+ik)/m)
   &=& \frac{\pi i}{2m}\sum_{k=1}^{m-1} g(k)\bigg(1+
   \frac{2}{e^{2(y+ik)\pi/m}-1}\bigg)\\
   &=& \tfrac12 \pi i M(g) + \frac{\pi i}{m}\sum_{k=1}^{m-1}
   \frac{g(k)}{e^{2\pi(y+ik)/m}-1}\\
   &=& \frac{\pi i}{m}\sum_{k=1}^{m-1}
   \frac{g(k)}{e^{2\pi y/m}\w^k-1},
\end{eqnarray*}
which proves~(\ref{ToddOmega}).

Finally, define $E:\C\to\C$ by $E(z):=\exp(2\pi z/m)$.
Then~(\ref{ToddOmega}) can be restated as
\begin{eqnarray*}
   T_g(-1,y)
   &=& \frac{\pi i}{m}\sum_{k=1}^{m-1}\frac{g(k)}{E(y+ik)-1}\\
   &=& \frac{\pi i}{2m}
   \sum_{k=1}^{m-1}\bigg(\frac{g(k)}{E(y+ik)-1}
   +\frac{g(m-k)}{E(y+i(m-k))-1}\bigg)\\
   &=& \frac{\pi i}{2m}\sum_{k=1}^{m-1}\bigg(\frac{1}{E(y+ik)-1}
       -\frac{1}{E(y-ik)-1}\bigg)g(k)\\
   &=& \frac{\pi
   i}{2m}\sum_{k=1}^{m-1}\bigg(\frac{E(y-ik)-E(y+ik)}{E(2y)+1
   -E(y+ik)-E(y-ik)}\bigg)g(k)\\
   &=& \frac{\pi
   i}{2m}\sum_{k=1}^{m-1}\bigg(\frac{E(-ik)-E(ik)}
              {E(y)+E(-y)-E(ik)-E(-ik)}\bigg)g(k)\\
   &=& \frac{\pi}{2m}\sum_{k=1}^{m-1}\frac{g(k)\sin(2\pi k/m)}{
   \cosh(2\pi y/m)-\cos(2\pi k/m)},
\end{eqnarray*}
which is~(\ref{ToddReal}). \eop

\begin{Lem}\label{Lem:6}
Let $m$ be a positive integer, $\w:=\exp(2\pi i/m)$, $y$ a
positive real number, $g:\Z\to\C$ even and periodic of period $m$,
$M(g)$ as in~(\ref{MeanValue}) and $T_g$ as in~(\ref{Tgdef}).
Furthermore, let $A_k(y) := \exp(2\pi y/m)-\cos(2\pi k/m)$.  Then
\begin{eqnarray}
  T_g(0,y) &=& \sum_{k=1}^\infty \frac{g(k)}{k^2+y^2}\nonumber\\
   &=& -\tfrac12 y^{-2}g(0)+\tfrac12 \pi y^{-1}m^{-1}
      \sum_{k=0}^{m-1}g(k)\coth(\pi(y+ik)/m)\label{TevenCoth}\\
    &=& \tfrac12 \pi y^{-1} M(g)-\tfrac12 y^{-2}g(0)
     + \frac{\pi}{y m}\sum_{k=0}^{m-1} \frac{g(k)}{e^{2\pi
     y/m}\w^k-1}\label{TevenOmega}\\
   &=& \tfrac12 \pi y^{-1} M(g)-\tfrac12 y^{-2}g(0)
     - \frac{\pi}{y m}\sum_{k=0}^{m-1} \frac{g(k) A_k(-y)}
     {A_k(y)+A_k(-y)}.\label{TevenReal}
\end{eqnarray}
\end{Lem}

\noindent{\bf Proof.} Since $g$ is even,
\begin{eqnarray*}
   2y\sum_{k=1}^\infty\frac{g(k)}{k^2+y^2}
   &=& \sum_{k=1}^\infty
   \bigg(\frac{1}{y+ik}+\frac{1}{y-ik}\bigg)g(k)\\
   &=& \sum_{k=1}^\infty \bigg(\frac{g(k)}{y+ik}+
   \frac{g(-k)}{y-ik}\bigg)\\
   &=& P.V. \sum_{0\ne k\in\Z}\frac{g(k)}{y+ik}.
\end{eqnarray*}
Formula~(\ref{TevenCoth}) now follows directly from
Lemma~\ref{Lem:4}.

If we now write
\begin{eqnarray*}
    \sum_{k=0}^{m-1}g(k)\coth(\pi(y+ik)/m)
    &=& \sum_{k=0}^{m-1}
    g(k)\bigg(1+\frac{2}{e^{2\pi(y+ik)/m}-1}\bigg)\\
    &=& m M(g) + \sum_{k=0}^{m-1} \frac{2 g(k)}{e^{2\pi
    y/m}\w^k-1},
\end{eqnarray*}
we see that formula~(\ref{TevenOmega}) follows
from~(\ref{TevenCoth}).

Finally, let $E:\C\to\C$ be defined by $E(z):=\exp(2\pi z/m)$.
Then
\begin{eqnarray*}
   \sum_{k=0}^{m-1} \frac{g(k)}{e^{2\pi
    y/m}\w^k-1} &=& \frac12 \sum_{k=1}^{m-1}\bigg(\frac{g(k)}{E(y+ik)-1}
      + \frac{g(m-k)}{E(y+i(m-k))-1}\bigg)\\
   &=&\frac12 \sum_{k=1}^{m-1}\bigg(\frac{1}{E(y+ik)-1}
      +\frac{1}{E(y-ik)-1}\bigg)g(k)\\
   &=&\frac12 \sum_{k=1}^{m-1}\bigg(\frac{E(y-ik)+E(y+ik)-2}
      {E(2y)+1-E(y+ik)-E(y-ik)}\bigg)g(k)\\
   &=&\frac12\sum_{k=1}^{m-1}\bigg(\frac{E(ik)+E(-ik)-2E(-y)}
      {E(y)+E(-y)-E(ik)-E(-ik)}\bigg)g(k)\\
   &=&\frac12\sum_{k=0}^{m-1}\bigg(\frac{2\cos(2\pi
   k/m)-2\exp(-2\pi y/m)}{2\cosh(2\pi y/m)-2\cos(2\pi
   k/m)}\bigg)g(k)\\
   &=& -\sum_{k=0}^{m-1}\frac{g(k)A_k(-y)}{A_k(y)+A_k(-y)}
\end{eqnarray*}
shows that~(\ref{TevenReal}) follows from~(\ref{TevenOmega}).
\eop

\begin{Lem}\label{Lem:7} Let $m$ and $q$ be positive integers, $x$ a
positive real number, $g:\Z\to\C$ odd and periodic of period $m$,
and $L(s,g)$ as in~(\ref{Lsg}). Then
\begin{eqnarray*}
   &&\tfrac12\pi x L(2q,g) + \pi x\sum_{n=1}^\infty
   \frac{n^{-2q}g(n)}{e^{2\pi n/x}-1}
   + \sum_{j=0}^q (-1)^{j}x^{-2j+2}\z(2j)L(2q-2j+1,g)\\
   &=& (-1)^q x^{-2q+2}\pi i m^{-1}\sum_{k=1}^{m-1} g(k)
      \sum_{n=1}^\infty \frac{n^{-2q}}{e^{2\pi nx/m}\w^k-1}\\
   &=&
   \tfrac12(-1)^q x^{-2q+2}\pi m^{-1}\sum_{k=1}^{m-1}
      g(k)\sin(2\pi k/m)
    \sum_{n=1}^\infty\frac{n^{-2q}}{\cosh(2\pi nx/m)-\cos(2\pi
    k/m)}.
\end{eqnarray*}
\end{Lem}

\noindent{\bf Proof.}  In Lemma~\ref{Lem:3}, set $s=2q-1$ and let
$f$ be the restriction of $g$ to the set $\Z^+$ of positive
integers, so that $L(s,g)=F(s)$ and $T_g=T_f$.   As in the proof
of Lemma~\ref{Lem:5}, since $g$ is odd and has period $m$,
$M(g)=0$, and thus $\sigg=\sigf\le 0$. Therefore, the convergence
condition $\Re(s)>\sigg+2q-2$ of Lemma~\ref{Lem:3} is satisfied.
Substituting the formulas of Lemma~\ref{Lem:5} for $T_g(-1,nx)$
completes the proof. \eop

\begin{Lem}\label{Lem:8}
Let $m$ and $q$ be positive integers, $x$ a positive real number,
$g:\Z\to\C$ even and periodic of period $m$, $M(g)$ as
in~(\ref{MeanValue}), and $L(s,g)$ as in~(\ref{Lsg}). Then
\begin{eqnarray*}
   &&\tfrac 12\pi x L(2q+1,g)+\pi x\sum_{n=1}^\infty
   \frac{n^{-2q-1}g(n)}{e^{2\pi n/x}-1}
    +\sum_{j=0}^{q+1} (-1)^{j}x^{-2j+2}\z(2j)L(2q-2j+2,g)\\
    &=& (-1)^q  \pi x^{-2q+1} \bigg\{ \tfrac12
    M(g)\z(2q+1) + \frac1m \sum_{k=0}^{m-1} g(k)
    \sum_{n=1}^\infty \frac{n^{-2q-1}}{e^{2\pi
    nx/m}\w^k-1}\bigg\}\\
    &=& \tfrac12(-1)^q\pi x^{-2q+1} \bigg\{  M(g)\z(2q+1)\\
       &&\qquad  - \frac1m\sum_{k=0}^{m-1} g(k)
   \sum_{n=1}^\infty \bigg(\frac{\exp(-2\pi nx/m)-\cos(2\pi k/m)}
   {\cosh(2\pi nx/m)-\cos(2\pi k/m)}\bigg)n^{-2q-1}\bigg\}.
\end{eqnarray*}
\end{Lem}

\noindent{\bf Proof.}  In Lemma~\ref{Lem:3}, set $s=2q$ and let
$f$ be the restriction of $g$ to set $\Z^+$ of positive integers.
Since $\sigg\le 1$, the convergence condition $\Re(s)>\sigg+2q-2$
of Lemma~\ref{Lem:3} is satisfied.  Next, substitute the formulas
of Lemma~\ref{Lem:6} for $T_g(0,nx)$.  The term
\[
   -\tfrac12 (-1)^q x^{-2q} g(0)\z(2q+2)
\]
arising from the term $-g(0)/2y^2$ in~(\ref{TevenOmega})
and~(\ref{TevenReal}) can be absorbed into the sum on $j$, since
for even $g$, $L(0,g)=- g(0)/2$ by
Proposition~\ref{Prop:Lclosedcplx}. \eop

\noindent{\bf Proof of Theorems~\ref{Thm:gcplx}
and~\ref{Thm:greal}.}   Set $x=\pi/\al$, $\be=\pi x=\pi^2/\al$ in
Lemma~\ref{Lem:7} and Lemma~\ref{Lem:8}.  \eop

\end{document}